# $K$-Exit for OD-Characterization of the finite Groups

Majid Akbari

*Department of Mathematics, Payame Noor University,*
*Tehran, Iran,*
*E-mail*: `majid.math@pnu.ac.ir`

June 6, 2025

## Abstract

Suppose $G$ is a finite group, and $K$ is a normal solvable subgroup of $G$. In this paper, we present a powerful method which is called $K$-Exit. Using this method, we find some members of $\pi(G)$ that are not members of $\pi(K)$ and for this purpose, we introduce the tables of $K$-Exit. This method can be applied to OD-Characterization of the finite Groups.

## 1 Introduction

The prime graph associated with the finite group $G$ is a graph with a set of vertices $\pi(G)$ and two vertices $p$ and $q$ are joined if and only if group $G$ has an element of order $pq$. We should mention that the degrees of the vertices of the prime graph are in ascending order and they are known as the degree pattern of $G$.

Suppose that $M$ is a finite group. If each finite group with the same order and the same degree pattern as $M$ is isomorphic to $M$, then we say that $M$ is OD-characterization [1]. For example, all simple finite groups of Lie type $L_2(q)$ where $q \neq 2, 3$ are OD-characterization [1, 8] and all alternating simple finite groups $A_n$ where $n = p, p+1, p+2$, and $p$ is a prime number greater than or equal to 5 are OD-characterization [1, 2], too. Also, all alternating simple finite groups $A_n$ that $5 \leq n \leq 100, n \neq 10$ are OD-characterization [3, 4, 5, 6, 7].

In many papers, to show the OD-characterization of finite group $M$, it is necessary to find members of $\pi(G)$ which are not belong to $\pi(K)$, where $G$ is a finite group with the same order and the same degree pattern as $M$ and $K$ is the maximal normal solvable subgroup of $G$. If we find more and more members of $\pi(G)$ which do not belong to $\pi(K)$, then we have more chance to showing that $M$ is OD-characterization. For this purpose, in the following, we present the $K$-Exit method and introduce the tables of $K$-Exit.

## 2 lemmas and theorems

In this section, we prove the lemmas and theorems needed to present the tables of K-Exit.

---

**Lemma 2.1** *Let $G$ be a finite group. Let $p, q \in \pi(G)$ and $|G| = p^i q$. If $q \nmid p^i - 1$ then $G$ has a subgroup of order $pq$.*

*Proof.* Suppose $P \in Syl_p(G), Q \in Syl_q(G)$. We show that $Q \lneq N_G(Q)$. If $Q = N_G(Q)$, then $|Syl_q(G)| = [G : Q] = p^i$. Now, by sylow theorems $q \mid p^i - 1$, which is a contradiction. Thus $Q \lneq N_G(Q)$, and we conclude that $|Q| \lneq |N_G(Q)|$. Therefore $p \mid |N_G(Q)|$. Thus, there is $t \in N_G(Q)$ such that $o(t) = p$. Now, it is clear that $Q\langle t \rangle$ is a group of order $pq$. and the proof is complete. □

**Definition 2.2** *Let $m \in N$. We denote by $[m]$ the least common multiple of $1, 2, ..., m-1, m$. In other words, $[m] = [1, 2, ..., m-1, m]$.*

**Lemma 2.3** *Let $G$ be a finite group and $p, q \in \pi(G)$. suppose that $m \in N$ and there exists $1 \leqslant i \leqslant m$ such that $|G| = p^i q$. If $q \nmid p^{([m], q-1)} - 1$, then $G$ has a subgroup of order $pq$.*

*Proof.* $p^{([m], q-1)} - 1 = (p^{[m]} - 1, p^{q-1} - 1)$. Thus $q \nmid p^{[m]} - 1$ or $q \nmid p^{q-1} - 1$. According to the Fermat's theorem $q \mid p^{q-1} - 1$. Thus $q \nmid p^{[m]} - 1$, and we conclude that $q \nmid p^i - 1$. Now, by lemma 2.1, $G$ has a subgroup of order $pq$. □

**Theorem 2.4** *Let $G$ be a finite group and $p, q \in \pi(G)$. suppose that $m, n \in N$ and there exists $1 \leqslant i \leqslant m$ such that $|G| = p^i q^n$. If $p \nmid q^n - 1$, $q \nmid p^{([m], q-1)} - 1$ then $G$ has an abelian subgroup of order $pq$.*

*Proof.* Suppose $P \in Syl_p(G), Q \in Syl_q(G)$. We show that $P \lneq N_G(P)$. If $P = N_G(P)$, then $|Syl_p(G)| = [G : P] = q^n$. Now, by sylow theorems $p \mid q^n - 1$, which is a contradiction. Thus $P \lneq N_G(P)$, and we conclude that $|P| \lneq |N_G(P)|$. Therefore $q \mid |N_G(P)|$. So, there is $t \in N_G(P)$ such that $o(t) = q$. It is clear that $P\langle t \rangle$ is a group of order $p^i q$. Now, since $q \nmid p^{([m], q-1)} - 1$, by the lemma 2.3, $G$ has a subgroup $H$ of order $pq$. On the other hand $p \nmid q - 1$, $q \nmid p - 1$. Thus the group $H$ is abelian, and the proof is complete. □

**Definition 2.5** *Suppose that $G$ is a finite group and $p \in \pi(G)$. We denote the Power of $p$ in $G$ by $w_G(p)$ and introduce it as follow:*

$$w_G(p) = max\{i \mid p^i \mid |G|\}$$

**Definition 2.6** *Suppose that $G$ is a finite group and $p \in \pi(G)$. We introduce $L(p, G)$ as follow:*
$$L(p, G) = \{q \in \pi(G) \setminus \{p\} \mid p \nmid q^n - 1, q \nmid p^{([m], q-1)} - 1, where \ \ m = w_G(p), n = w_G(q)\}$$

**Theorem 2.7** *Suppose that $G$ is a finit group and $K$ is a normal solvable subgroup of $G$ and $GK(G)$ is the prime graph of $G$. Let $p \in \pi(G)$. If $d_G(p) < |L(p, G)|$, then $p \notin \pi(K)$.*

*Proof.* Suppose this proposition is fals. In other words, $p \in \pi(K)$. Let $q \in L(p, G)$ and $w_G(p) = m, w_G(q) = n$. Let $w_K(p) = i$. It is clear that $1 \leqslant i \leqslant m$. Assume $Q \in Syl_q(G)$. $K$ is a normal solvable and $Q$ is solvable. Thus $KQ$ is solvable, and we conclude that $KQ$ has a Hall $\{p, q\}$-subgroup $F$. Clearly, $|F| = p^i q^n$. On the other hand $q \in L(p, G)$. Hence $p \nmid q^n - 1, q \nmid p^{([m], q-1)} - 1$. Now by theorem 2.4, $F$ has a abelian subgroup of order $pq$. Therefore $p \sim q$, and we conclude that $p$ is connected to each vertex in $L(p, G)$. So $d_G(p) \geqslant |L(p, G)|$, which is a contradiction. Thus $p \notin \pi(K)$. □

**Remark 2.8**

Suppose that $G$ is a finit group and $K$ is a normal solvable subgroup of $G$ and $GK(G)$ is the prime graph of $G$. Theorem 2.4 is useful for checking that the two vertices in $\pi(K)$ are joined. $[m]$ causes the calculations to be reduced to one step. On the other hand $[m]$ is large. Therefore calculations will be long and tedious, but the $([m], q-1)$ is relatively small and decreases the calculations. The theorem 2.7 is useful for checking whether a specific members of $\pi(G)$ is outside of $\pi(K)$ or not. In the following, we will present other theorems to check whether the members of $\pi(G)$ are outside of $\pi(K)$ or not. Also, we introduce the tables of $K$-Exit.

**Theorem 2.9** *Let $G$ be a finite group and $p, q \in \pi(G)$. Let $m, n \in N$ and suppose that there exists $1 \leqslant i \leqslant m$ such that $|G| = p^i q^n$. If $p \nmid q^n - 1$, $q \nmid p^j - 1$ for all $1 \leqslant j \leqslant m$, then $G$ has an abelian subgroup of order $pq$.*

*Proof.* Suppose $P \in Syl_p(G), Q \in Syl_q(G)$. We show that $P \lneqq N_G(P)$. If $P = N_G(P)$, then $|Syl_p(G)| = [G : P] = q^n$. Now, by sylow theorems $p \mid q^n - 1$, which is a contradiction. Thus $P \lneqq N_G(P)$, and we conclude that $|P| \lneqq |N_G(P)|$. Therefore $q \mid |N_G(P)|$. Thus, there is $t \in N_G(P)$ such that $o(t) = q$. It is clear that $P\langle t \rangle$ is a group of order $p^i q$. Now, since $q \nmid p^i - 1$, by lemma 2.1, $G$ has a subgroup $H$ of order $pq$. On the other hand $p \nmid q - 1$, $q \nmid p - 1$, Thus the group $H$ is abelian, and the proof is complete. $\square$

**Definition 2.10** Suppose that $G$ is a finite group and $p \in \pi(G)$, $w_G(p) = m$. We introduce $\theta(p)$, $\overline{\theta}(p)$ as follow:

$$\theta(p) = \pi(G) \setminus \{p\} \cup (\bigcup_{i=1}^{m} \pi(p^i - 1))$$

$$\overline{\theta}(p) = \pi(G) \setminus \{p\} \cup \pi(p^m - 1))$$

**Definition 2.11** Suppose that $G$ is a finite group and $p \in \pi(G)$. We denote the page $p$ in $G$ by $H(p, G)$ and introduce it as follow:

$$H(p, G) = \{q \in \theta(p) \mid p \in \overline{\theta}(q)\}$$

**Theorem 2.12** *Suppose that $G$ is a finit group and $K$ is a normal solvable subgroup of $G$ and $GK(G)$ is the prime graph of $G$. Let $p \in \pi(G)$. If $d_G(p) < |H(p, G)|$, then $p \notin \pi(K)$.*

*Proof.* Suppose this proposition is fals. In the other words, $p \in \pi(K)$. Let $q \in H(p, G)$, $w_G(p) = m, w_G(q) = n$, and $w_K(p) = i$. It is clear that $1 \leqslant i \leqslant m$. Assume $Q \in Syl_q(G)$. $K$ is a normal solvable and $Q$ is solvable. Thus $KQ$ is solvable, and we conclude that $KQ$ has a Hall $\{p, q\}$-subgroup $F$. Clearly, $|F| = p^i q^n$. On the other hand $q \in H(p, G)$. Hence $p \in \overline{\theta}(q)$, $q \in \theta(p)$, and it follows that $p \nmid q^n - 1$, $q \nmid p^j - 1$, for all $1 \leqslant j \leqslant m$. Now by theorem 2.9, $F$ has a abelian subgroup of order $pq$. Therefore $p \sim q$, and we conclude that $p$ is connected to each vertex in $H(p, G)$. So $d_G(p) \geqslant |H(p, G)|$, which is a contradiction. Thus $p \notin \pi(K)$. $\square$

## 3 The tables of $K$-Exit and Examples

Suppose that $G$ is a finit group and $K$ is a normal solvable subgroup of $G$. Theorem 2.12 helps us to find some members of $\pi(G)$ that don't belong to $\pi(K)$ and helps us to design the tables of $K$-Exit. For this purpose, we arrange the members of $\pi(G)$ in the first column and assign the rows of each of them to $\theta$, $\bar{\theta}$, and $H$. Then we compare $|H|$ and the degree for each of them and use the Theorem 2.12.

**Example 3.1** $M = U_3(31)$.

OD-Characterization of the finite Groups $U_3(31)$ has been proven. for this purpose, it is necessary to find members of $\pi(G)$ that don't belong to $\pi(K)$ where $G$ is a finite group with the same order and degree pattern as $M$ and $K$ is the maximal normal solvable subgroup of $G$. See [9]. Here we do it by the tables of $K$-Exit.

$$|G| = |U_3(31)| = 2^{11} \cdot 3 \cdot 5 \cdot 7^2 \cdot 19 \cdot 31^3$$

$$D(G) = D(U_3(31)) = (3, 2, 2, 1, 1, 1),$$

| $p$ | $\theta(p)$ | $\bar{\theta}(p)$ | $H(p, G)$ | $d_G(p)$ | $|H(p, G)|$ | result |
|---|---|---|---|---|---|---|
| 2 | | | − | 3 | − | − |
| 3 | 5, 7, 19, 31 | 5, 7, 19, 31 | 5 | 2 | 1 | − |
| 5 | 3, 7, 19, 31 | 3, 7, 19, 31 | 3, 7, 19 | 2 | 3 | $5 \notin \pi(K)$ |
| 7 | 5, 19, 31 | 5, 19, 31 | 5, 19, 31 | 1 | 3 | $7 \notin \pi(K)$ |
| 19 | 5, 7, 31 | 5, 7, 31 | 5, 7, 31 | 1 | 3 | $19 \notin \pi(K)$ |
| 31 | 7, 19 | 7, 19 | 7, 19 | 1 | 2 | $31 \notin \pi(K)$ |

Table1. $K$-Exit for $U_3(31)$

**Example 3.2** $M = U_4(89)$.

OD-Characterization of the finite Groups $U_4(89)$ has been proven. See [9]. Suppose $G$ is a finite group with the same order and degree pattern as $M$, and $K$ is the maximal normal solvable subgroup of $G$. Now we use the $K$-Exit method.

$$|G| = |U_4(89)| = 2^9 \cdot 3^7 \cdot 5^3 \cdot 7 \cdot 11^2 \cdot 17 \cdot 89^6 \cdot 233 \cdot 373$$

$$D(G) = D(U_4(89)) = (6, 6, 6, 3, 6, 3, 4, 3, 3),$$

| $p$ | $\theta(p)$ | $\bar{\theta}(p)$ | $H(p,G)$ | $d_G(p)$ | $|H(p,G)|$ | result |
|---|---|---|---|---|---|---|
| 2 | | | − | 6 | − | − |
| 3 | 17, 89, 233, 373 | 5, 7, 11, 17, 89, 233, 373 | 17, 233 | 6 | 2 | − |
| 5 | 7, 11, 17, 89, 233, 373 | 3, 7, 11, 17, 89, 233, 373 | 7, 17, 89, 233, 373 | 6 | 5 | − |
| 7 | 5, 11, 17, 89, 233, 373 | 5, 11, 17, 89, 233, 373 | 5, 11, 17, 233, 373 | 3 | 5 | $7 \notin \pi(K)$ |
| 11 | 7, 17, 89, 233, 373 | 7, 17, 89, 233, 373 | 7, 17, 233, 373 | 6 | 4 | − |
| 17 | 3, 5, 7, 11, 89, 233, 373 | 3, 5, 7, 11, 89, 233, 373 | 3, 5, 7, 11, 89, 233, 373 | 3 | 7 | $17 \notin \pi(K)$ |
| 89 | − | 5, 7, 233 | − | 4 | − | − |
| 233 | 3, 5, 7, 11, 17, 89, 373 | 3, 5, 7, 11, 17, 89, 373 | 3, 5, 7, 11, 17, 89, 373 | 3 | 7 | $233 \notin \pi(K)$ |
| 373 | 5, 7, 11, 17, 89, 233 | 5, 7, 11, 17, 89, 233 | 5, 7, 11, 17, 233 | 3 | 5 | $373 \notin \pi(K)$ |

Table2. $K$-Exit for $U_4(89)$